%% file: ifac2010.tex
\documentclass{ifacconf}
\usepackage{amsfonts}
\usepackage{amsmath}
\usepackage{amssymb}
\usepackage{graphicx}      % include this line if your document contains figures
\renewcommand{\e}{\varepsilon}
\renewcommand{\d}{\delta}
\newcommand{\pd}[1]{\frac{\partial}{\partial #1}}
\newcommand{\RR}{\mathbb{R}}
\renewcommand{\SS}{\mathbb{S}}
\newcommand{\G}{\textrm{G}}
\newcommand{\E}{\cal{E}}
\newcommand{\N}{\cal{N}}

\newcommand{\f}{\textrm{f}}
\newcommand{\tf}{t_\textrm{f}}

\newcommand{\xf}{x_\f}

\newcommand{\ti}{t_0}
\newcommand{\intr}{\textrm{int}}
\newcommand{\vol}{\textrm{vol}}
\DeclareMathOperator{\maximize}{\textrm{max.}}
\DeclareMathOperator{\minimize}{\textrm{min.}}

\newcommand{\subjto}{\textrm{subj. to}}

 \newcommand{\goal}{{\cal G}}
\newcommand{\fun}{{\cal F}}

\begin{document}
\begin{frontmatter}

  \title{Invariant Funnels around Trajectories using Sum-of-Squares Programming. \thanksref{footnoteinfo}} 

\thanks[footnoteinfo]{This material is based upon work supported by the National Science Foundation under Grant No. 0915148.}

\author[CSAIL]{Mark M. Tobenkin}
\author[CSAIL]{Ian R. Manchester}
\author[CSAIL]{Russ Tedrake}
\address[CSAIL]{Computer Science and Artificial
  Intelligence Laboratory, Massachusetts Institute of Technology,
  Cambridge, MA 02129 USA (e-mail: mmt@mit.edu).}

\begin{abstract}                % Abstract of not more th an 250 words.
  This paper presents numerical methods for computing regions of
  finite-time invariance (funnels) around solutions of polynomial
  differential equations.  First, we present a method which exactly
  certifies sufficient conditions for invariance despite relying on
  approximate trajectories from numerical integration.  Our second
  method relaxes the constraints of the first by sampling in time. In
  applications, this can recover almost identical funnels but is much faster to
  compute.  In both cases, funnels are verified using Sum-of-Squares
  programming to search over a family of time-varying polynomial
  Lyapunov functions. Initial candidate Lyapunov functions are
  constructed using the linearization about the trajectory, and
  associated time-varying Lyapunov and Riccati differential
  equations. The methods are compared on stabilized trajectories of a
  six-state model of a satellite.
\end{abstract}


\end{frontmatter}
%===============================================================================

%\section{Introduction}
%\include{ifac2010intro.tex}
\section{Introduction}
\input{ifac2010intro}
\section{Preliminaries}
\input{ifac2010posinv}
\section{Optimization Procedure}
\input{ifac2010sos}
\section{Sampling in Time}
\input{ifac2010samp}

\section{Numerical Experiments}
\input{ifac2010numerics}
\section{Discussion and Future Work}
\input{ifac2010discussion}
\bibliography{elib}
\end{document}

%% file: ifac2010intro.tex
In this paper we propose algorithms for computing a finite-time ``funnel'' of a dynamical system: a set of initial conditions whose solutions are guaranteed to enter a certain goal region at a particular time. This work is motivated by new algorithms for nonlinear control design wherein the controller is constructed from a tree of finite-time trajectories, each locally stabilized and all leading to a certain goal point (see \cite{Tedrake10}). Sparseness of this tree is advantageous, and is directly related to the size of the funnel that can be verified for each trajectory in the tree.  % The authors have also adapted the
% procedure proposed in this paper to verify regions of attraction of
% limit cycles (see \cite{Manchester10a, Manchester10b}).

The basic method is to search over a class of time-varying Lyapunov functions about the trajectory, and verify a positive-invariance condition via Sum-of-Squares programming. A preliminary algorithm for this verification was suggested in \cite{Tedrake10}. This paper extends that work by replacing a greedy set of nonconvex optimizations by an alternation of convex optimizations (similar to that proposed for equilibria in \cite{Jarvis-Wloszek05}).  It is shown that the verified regions can be made exact even if the trajectory is an approximate (numerical) solution. Furthermore, an alternative time-sampled verification is suggested which recovers almost identical funnels and is substantially faster to compute.

Recently, a great deal of research has been dedicated to
region-of-attraction analysis on polynomial vector fields
(\cite{Papachristodoulou02, Parrilo03a, Topcu08, Tan08}) and more
general non-polynomial vector fields
(\cite{Papachristodoulou05,Chesi09}) through Sum-of-Squares
programming. There is comparatively little written in the
Sum-of-Squares literature addressing time-varying systems. In
\cite{Julius09} finite time-invariance around trajectories is explored
to provide outer approximations of the set reached from some initial
conditions.  This is accomplished by using regionally valid Lyapunov
certificates to construct barrier functions bounding exponentia, certifying solutions do
not enter keep-out sets.  By contrast, the algorithm in this paper
constructs inner-approximations of the solutions which can reach a
goal region through the construction of time-varying Lyapunov
functions.

\subsection{Notation}
% Given a vector $x \in \RR^n$, $x'$ denotes transpose.
% The set of symmetric positive definite matricies in $\RR^{n \times n}$
% is denoted $\SS_+^n$.  Given $P \in \SS_+^n$ we define $\|x\|_P^2 :=
% x'Px$. For $X,Y \in \RR^{n\times n}$, $X > Y$ (resp. $X \geq Y$)  indidcates $X-Y$ is
% positive definite (resp. positive semidefinite).
We denote the set of $n$-by-$n$ positive definite matrices by
$\SS_+^n$.  For $P \in \SS_+^n$ and $x \in \RR^n$ we use $\|x\|_{P}^2$
as short hand for $x'Px$.  We denote the set of polynomials in $x \in
\RR^n$ with real coefficients by $\RR[x]$.  The subset of these
polynomials which are Sum-of-Squares (SOS) is denoted $\Sigma[x]$,
that is: $p(x) \in \Sigma[x]$ if and only if there exists
$\{g_i(x)\}_{i=1}^k \subset \RR[x]$ such that $p(x) = \sum_{i=1}^k
g_i(x)^2$.  We similarly denote polynomials and SOS polynomials in
multiple vector valued variables $(x_1,\ldots,x_k) \in \RR^{n_1}\times
\ldots \times \RR^{n_k}$ by $\RR[x_1,\ldots,x_k]$ and
$\Sigma[x_1,\ldots,x_k]$.  Finally, for a set $A \subset \RR^n$,
$\intr(A)$ refers to its interior.

% When referring to functions $p(t,x)$ defined on $[a,b]\times \RR^n$
% as piecewise polynomial we will commonly denote the ``knots'' of the
% function as $\{t_i\}_{i=0}^N$ with $a = t_0 < t_1 < \ldots <
% t_N = b$, and the ``pieces'' as $p_i(t,x) \in \RR[t,x]$. The
% function is then defined by:
% $$p(t,x) = \begin{cases} p_i(t,x) & t \in (t_{i-1},t_{i}] \\ p_1(t,x)
% & t=t_0 \end{cases}.$$
% The function is clearly continuous in $t$ if for $i = \{1,\ldots,N-1\}$,
% $p_i(t_i,x) = p_{i+1}(t_{i},x)$.  

%%% Local Variables: 
%%% mode: latex
%%% TeX-master: "ifac2010"
%%% End: 

%% file: ifac2010posinv.tex
\label{sec:posinv}
We are given a time-varying nonlinear dynamical system:
\begin{equation}\label{eq:dyn}
\frac{d}{dt} x(t) = f(t,x(t)),
\end{equation}
where $f: [\ti,\tf] \times \RR^n \mapsto \RR^n$ is piecewise
polynomial in $t$ and polynomial in $x$.  Further we are given an
 bounded ``goal region'', $\goal \subset \RR^n$, with non-empty interior.  In applications, we will generally
be investigating systems \eqref{eq:dyn} arising as closed loop systems
stabilizing some trajectory.   For convenience
we make the following definitions. 

\begin{defn}
  Given the dynamics \eqref{eq:dyn}, a set $\fun \subset [\ti,\tf] \times \RR^n$ is a {\it funnel} if
  for each   $(\tau, x_\tau)$ in $\fun$, the solution to
  \eqref{eq:dyn} with $x(\tau) =x_\tau$ exists on $[\tau,\tf]$ and for
  each $t \in [\tau,\tf]$ we have $(t,x(t)) \in \fun$.
\end{defn}
\begin{defn}
  Given the dynamics \eqref{eq:dyn} and the goal region $\goal$, a set $\fun \subset [\ti,\tf]
  \times \RR^n$ is a {\it funnel into $\goal$} if it is a funnel and for any $x \in \RR^n$
  we have $(\tf,x) \in \fun$ implying $x \in \goal$.
\end{defn}

We see that a funnel into $\goal$ is an inner-approximation of the set of
solutions which flow through $\goal$ at the time $\tf$.  In this work
we are interested in finding the largest possible funnel into $\goal$,
measured, for example, by volume as a subset in
$[\ti,\tf]\times\RR^n$.  Our approach uses time-varying Lyapunov
functions, exploiting the following Lemma.

% Here $u: [\ti,\tf] \mapsto
% \RR^m$ is drawn from a class $\U$ of piecewise continuous signals.
% In this paper, we seek to address specific instances of the following
% general problem. 

% \begin{prob}
%   Given a nonlinear dynamical system \eqref{eq:ctldyn}, and a bounded,
%   open set $\goal \subset \RR^n$, let $\fun^\star \subset [\ti,\tf]
%   \times \RR^n$ be defined by  $(\tau,x_\tau) \in
%   \fun^\star$ if and only if there exists a $u^\star \in \U$ such that
%   the solution $x(t)$ with $x(\tau) = x_\tau$ and $u(t) \equiv
%   u^\star(t)$ exists on $[\tau,\tf]$ and $x(\tf)$ lies inside $\goal$.
%   Find an inner approximation $\fun \subset \fun^\star$ which
%   is measureable, and has as large a volume in $[\ti,\tf]\times \RR^n$
%   as possible.
% \end{prob}

\begin{lem}\label{lem:posinv}
  Let $V(t,x)$ be a function $V: [\ti,\tf] \times \RR^n
  \mapsto [0,\infty)$ piecewise continuously differentiable\footnote{We
will generally write conditions on $\frac{\partial}{\partial t} V$.
At points of discontinuity, we require these conditions to hold for
both the left and right derivative of $V$ with respect to $t$.} with respect to $t$ and
  continuously differentiable with respect to  $x$.  For each $t \in [\ti,\tf]$ we define:
$$\Omega_{t} := \{ x \; |
\; V(t,x) \leq 1 \},$$
and
$$\partial \Omega_{t} := \{x \; | \; V(t,x) = 1\}.$$

% Given a function $\pi: [\ti,\tf] \times \RR^n \mapsto \RR^n$,
% continuous in $x$ and piecewise continuous in $t$, 
If for each $t \in [\ti,\tf]$ and $x \in \partial\Omega_{t}$ we have:
\begin{equation}\label{eq:posinv}
\frac{\partial}{\partial x}V(t,x) f(t,x) + \frac{\partial}{\partial
t}V(t,x) < 0, 
\end{equation}
then the set:
\begin{equation} \label{eq:F}
  \fun = \{ (t,x) \, | x \in \Omega_{t} \}
\end{equation}
is a funnel.  If additionally $\Omega_{\tf} \subset \goal$, then
$\fun$ is a funnel into $\goal$.
\end{lem}

% \begin{defn}
% Given a dynamical system \eqref{eq:dyn}, a set $\fun \subset
% [\ti,\tf] \times \RR^n$ will be called a {\it positive funnel} if 
% for each $(\tau,x_\tau) \in \fun$ there the
% solution $x(t)$ to \eqref{eq:dyn} with $x(\tau) = x_\tau$ exists on
% $[\tau,\tf]$ and the set of points $(t,x(t))$ lie inside $\fun$ for
% all $t \in [\tau,\tf]$.
% \end{defn}

% \begin{defn}
% Given a dynamical system \eqref{eq:dyn}, a positive funnel $\fun$
% and a set $\goal \subset
% \RR^n$, if $(\tf,x) \in \fun$ implies $x \in \goal$ we say that $\fun$
% is a {\it positive funnel into $\goal$}.
% \end{defn}

% This can be seen as a slight generalization of positive invariance to
% the extended state-space, $[\ti,\tf] \times \RR^{n}$.  It is useful to
% note that if $\fun_1$ and $\fun_2$ are positive funnels into $\goal$,
% clearly $\fun_1 \cup \fun_2$ is also a positive funnel.

To leverage the SOS relaxation, %  we now focus our attention to systems
% where $f(t,x,u)$ is polynomial in $x$ and $u$ and piecewise
% polynomial in $t$.  We similarly restrict ourselves to function
% $\pi(t,x)$ which are piecewise polynomial in $t$ and polynomial in
% $x$.  Finally,
our goal region $\goal$ will be the closed interior of an ellipse:
$$\goal = \{ x \; | \; \|x\|_{P_\G}^2 \leq 1\},$$
where $P_\G \in \RR^{n\times n}$ is a symmetric positive definite
matrix.  More general goal regions can be defined as sub-level sets of
polynomials.

% We are given a nonlinear dynamical system:
% \begin{equation}\label{eq:polyadyn}
% \dot x = f(t,x),
% \end{equation}
% with $x \in \RR^n$ and $f: [\ti,\tf] \times \RR^n \mapsto \RR^n$ being a
% polynomial function in $x$ and piecewise polynomial in $t$.  Further
% we are given an ellipsoidal goal
% region, $\goal$ specified a symmetric positive definite matrix
% $P_\G \in \RR^{n\times n}$:
% $$\goal = \{ x \; | \; x^TP_\G x < 1\}.$$

\subsection{A Class of Lyapunov Functions}
We begin by describing a class of candidate Lyapunov
functions based on solutions of \eqref{eq:dyn} and related Lyapunov
differential equations.  To leverage SOS we approximate these
solutions by piecewise polynomials in time.  This approximation does
not render our certificates inexact.  Further, we describe how
sufficiently fine approximation of these solutions guarantees the
existence of a certificate.

We assume we have  access to  a nominal solution
of \eqref{eq:dyn},
$x_0: [\tau,\tf] \mapsto \RR^n$ with $\tau \in [\ti,\tf)$,
% and nominal policy, $\pi_0: [\tau,\tf]
% \times \RR^n \mapsto \RR^m$,
such that $x_\f := x_0(\tf) \in \intr(\goal)$.  
We solve for a symmetric positive definite matrix $P_\f
\in \SS^n_+$ describing the largest volume ellipse centered on $\xf$,
$$\E_\f = \{ x \; | \; \|x-\xf\|_{P_\f}^2 \leq
1\},$$ contained in the goal region (i.e. $\E_\f \subset \goal$).
When $\goal$ is an ellipse, this containment problem is an SDP.  Given
a more general polynomial goal region, one can relax the problem to
finding the largest sphere centered on $\xf$ contained in
$\goal$ to a SOS program. 

Our Lyapunov functions are parameterized by a time-varying,
symmetric positive definite matrix $P^\star:~[\tau,\tf] \mapsto \SS^n_+$:
\begin{equation}\label{eq:idealV}
V^\star(t,x) = \| x - x_0(t) \|_{P^\star(t)}^2.
\end{equation}
We require $P^\star(\tf) \geq P_\f$,  so that $\Omega_{\tf}$, the one
sub-level set of $V^\star(\tf,x)$, is a subset of $\E_\f$ and thus contained in
the goal region.

Our optimization approach requires an initial feasible candidate
Lyapunov function which is then improved via an iterative process.  To
construct this initial candidate we look at the dynamics linearized
about the trajectory by constructing $A: [\tau,\tf] \mapsto
\RR^{n\times n}$, with $A(t) =  \frac{\partial f}{\partial
x}(t,x_0(t))$.  Then, fixing a function $Q: [\tau,\tf] \mapsto
\SS_+^n$, we solve the Lyapunov differential equation:
\begin{equation} \label{eq:lyap}
  -\dot P^\star_0(t) = A(t)' P^\star_0(t) + P^\star_0(t)A(t) + Q(t), \: P^\star_0(\tf) = P_\f
\end{equation}
over the interval $[\tau,\tf]$.  The following Lemma suggests a
procedure for finding a candidate Lyapunov function.

% \subsection{A Feasible Lyapunov Candidate}
% Our optimization approach requires an initial feasible candidate
% Lyapunov function which is then improved via an iterative process.  To
% construct this initial candidate we look at the dynamics linearized
% about the trajectory by constructing $A: [\ti,\tf] \mapsto
% \RR^{n\times n}$:
% \begin{equation}
%   A(\tau) = \left[\pd{x} f+
%   \pd{u}f\pd{x} \pi\right]_{t=\tau,x=x_0(\tau)}.
% \end{equation}
% and then, fixing some symmetric,  time-varying, positive definite
% matrix $Q: [\tau,\tf] \mapsto 
% \RR^{n\times n}$, solve the Lyapunov equation:
% \begin{equation} \label{eq:lyap}
%   -\dot P^\star_0(t) = A(t)^T P^\star_0(t) + P^\star_0(t)A(t) + Q(t), \: P^\star_0(\tf) = P_\f
% \end{equation}
% over the interval $[\tau,\tf]$.  The following Lemma suggests a
% procedure for finding a candidate Lyapunov function.

\begin{lem}
  Given a solution $P^\star_0(t)$ to the Lyapunov equation above, there
  exists a positive constant $c$ such that $V^\star$ defined in
  equation \eqref{eq:idealV} with:
  $$P^\star(t) = \exp\left(c\frac{\tf-t}{\tf-\tau}\right)P^\star_0(t)$$
  satisfies the conditions of Lemma \eqref{lem:posinv}, so that:
  $$\fun = \{(t,x)\in[\tau,\tf]\times\RR^n \; |\; V(t,x) \leq 1\}$$
  is a positive funnel into $\goal$.
\end{lem}

\begin{pf}
  We begin by changing coordinates to $\bar x(t) = x(t) - x_0(t)$, and
  defining $\bar f(t,\bar x)$ by:
  $$\dot {\bar x} = \bar f(t,\bar x) =
  f(t,x(t)) - f(t,x_0(t))$$
We can decompose $\bar f(t,\bar x)$ as:
\[
\dot x = A(t)\bar x+\tilde f(t,\bar x)
\]
where $\tilde f(t,\bar x)$ consists of the second and higher-order
terms in $\bar f(t,\bar x)$.
Taking 
\[
V(t,x) = \exp\left(c\frac{\tf-t}{\tf-\tau}\right)\bar x(t)'P^\star_0\bar x(t)
\]
we have
\begin{eqnarray}
\dot V(t,x) &=& \exp\left(c\frac{\tf-t}{\tf-\tau}\right)[2\bar x'P^\star_0\bar
f(t,\bar x)+\bar x'(\dot P^\star_0-cP^\star_0)\bar x]\notag\\
&=&\exp\left(c\frac{\tf-t}{\tf-\tau}\right)[2\bar x'P^\star_0\tilde f(t,\bar
x)-\bar x'(Q+cP^\star_0)\bar x].\notag
\end{eqnarray}
Now, $\partial \Omega_{t_f}$ is a compact set, so $\tilde f(t_f,x)$ is bounded for $x\in\partial \Omega_{t_f}$, therefore there exists a sufficiently large $c=c_1$ such that $\dot V(t_f ,x)<0$ for  $x\in\partial \Omega_{t_f}$. 

Since $P(t)>0$ and is continuous in $t$ and $\dot V$ is continuous in $x$ and $t$, there exists a time $t_m<t_f$ such that $\dot V(t,x)<0$ for all  $x\in\partial \Omega_{t}$ for all $t\in[t_m, t_f]$. 

We now show that this is also true on $t\in[\tau, t_f]$. Since
$t_m<t_f$, for any $\epsilon>0$ there exists a $c$ sufficiently large
that the sets $\partial \Omega_{t}$ are contained in the ball $|\bar x|<\epsilon$ for all $t\in[\tau ,t_m]$. 

Since $\bar x'P^*_0\tilde f(t,\bar x)$ contains third and higher
orders in $\bar x$, there exists a sufficiently small $\epsilon$ such
that $|\bar x(t)'P_0^\star(t)\tilde f(x,t)|<\bar
x(t)'(Q(t)+cP^\star_0(t))\bar x(t)$ for all $\bar x\in \partial \Omega_{t}$ for all $t\in[\tau, t_m]$. 

This implies that there exists a sufficiently large $c=c_2$  such that $\dot V(t,x)<0$ for all $x\in \partial \Omega_{t}$ for all $t\in[\tau, t_m]$.

Taking $c=\max\{c_1,c_2\}$ proves the Lemma.
\hfill$\Box$  
\end{pf}

As we only consider a finite time interval, this result guarantees
the funnel will have a non-empty intersection with $\{t\} \times
\RR^n$ for each $t \in [\tau,\tf]$.

\subsection{Polynomial Lyapunov Functions}
To exploit SOS programming, we develop an alternative Lyapunov
function to \eqref{eq:idealV} defined by a piecewise polynomial
functions: $\hat x_0: [\tau,\tf] \mapsto \RR^n$, $P_0: [\tau,\tf]
\mapsto \SS_+^n$ and $\rho: [\tau,\tf] \mapsto (0,\infty)$:
\begin{equation}\label{eq:V}
  V(t,x) = \|x-\hat x_0(t)\|_{P(t)}^2, \quad P(t) =
\frac{P_0(t)}{\rho(t)}.
\end{equation}
In particular, we have $\hat x_0(t)$ approximate $x_0(t)$ and $P_0(t)$
approximate $P_0^\star(t)$.  The function $\rho(t)$ is a time-varying
rescaling which we describe in the next section.

Note that the approximate nature of $\hat x_0(t)$ does not preclude
$V(t,x)$ from verifying a funnel as the conditions in Lemma
\ref{lem:posinv} are concerned with the {\it level-set} where
$V(t,x)=1$. As a result, the exact behavior of the function near
$V(t,x) = 0$ is not essential.

% A somewhat stronger statement is
% available.

% \begin{lem}
%   About uniform approximation of function and derivatives
%   leading to there always existing a funnel.
% \end{lem}
% \begin{pf} See Appendix.
% \end{pf}

% Again to exploit SOS programming, we approximate $P^\star(t)$.  In
% particular, we take:
% $$P(t) = \frac{P_0(t)}{\rho_0(t)}$$
% where $P_0(t)$ is a piecewise polynomial approximation of
% $P^\star_0(t)$ and $\rho_0: [\tau,\tf] \mapsto (0,\infty)$ is a
% piecewise polynomial in $t$ approximating:
% $$\rho_0(t) \approx \exp\left(-c\frac{\tf-t}{\tf-\tau}\right).$$
% The constraint that $P_0(t)$ be
% positive definite can ensured via a SOS program.  Similarly a SOS
% constraint can ensure the positivity of $\rho_0(t)$.

% \begin{equation}
%   A(\tau) = \left[\pd{x} f+
%   \pd{u}f\pd{x} \pi\right]_{t=\tau,x=x_0(\tau)}.
% \end{equation}
% To leverage the SOS decomposition, we approximate $x_0(t)$ by a
% continuous, piecewise polynomial function $\hat x_0(t)$ such that
% $\hat x_0(\tf) = x_0(\tf)$.  We then take:
% \begin{equation}\label{eq:V}
% V(t,x) = \| x - \hat x_0(t) \|_{P(t)}^2
% \end{equation}
% as our candidate Lyapunov function.

% \begin{rem}
%   Note that while $\hat x_0(t)$ is only an approximation of the
%   nominal solution, the conditions in Lemma \ref{lem:posinv} are
%   concerned with non-infintesimal {\it level-sets} of $V(t,x)$.  As a result,
%   exact behavior of the system near $V(t,x) = 0$ is not important.
% \end{rem}

%%% Local Variables: 
%%% mode: latex
%%% TeX-master: "ifac2010"
%%% End: 

%% file: ifac2010sos.tex
\label{sec:sos}
Once we restrict ourselves to a class of piecewise polynomial Lyapunov
functions, we can approach our optimization task as a bilinear
Sum-of-Squares program.  In particular, we will see that the
conditions of Lemma \eqref{eq:posinv} will be tests of polynomial
negativity on semi-algebraic sets.  To verify conditions on these
specific time-intervals and level-sets, we make use of the polynomial ${\cal
S}$-procedure (see \cite{Parrilo03a}).  We arrive in a program with
constraints bilinear in our parameterization of $V(t,x)$ and 
multipliers used in the ${\cal S}$-procedure.  In particular, for this work
we parameterize $V$ solely by our choice of time-varying $\rho(t)$.

Let $\{t_i\}_{i=0}^N$ be a set of knot points which contains the knot
points of $f(t,x)$, $\hat x_0(t)$, and $P_0(t)$.  In particular, we order the
knots so that $t_i < t_{i+1}$ for $i \in \N = \{0,\ldots,N-1\}$.

We define:
\begin{equation}
  \bar V(t,x) = \| x - \hat x_0(t)\|_{P_0(t)}^2 
\end{equation}
For each interval $[t_i,t_{i+1}]$ we define Lagrange multipliers
$\ell_i, \mu_i \in \RR[t,x]$.  Let $\rho_i$, $f_i$ and $\bar V_{i}$ be the
polynomial pieces of $\rho$, $f$ and $\bar V$ on the these intervals. 
For some constant $\e > 0$,
we will attempt to optimize a cost function $h(\rho)$ according to:
\begin{flalign}\label{eq:bilinear}
  \maximize_{\{\rho_i,\ell_i,\mu_i\}_{i\in\N}} \quad & h(\rho)\\
  \subjto \quad & \rho_{N-1}(\tf) = 1,
  \nonumber\\
  &\forall \; i \in \N:\nonumber\\
  & \rho_i(t) \in \Sigma[t], \rho_i(t_{i+1}) = \rho_{i+1}(t_{i+1}),\nonumber\\
  & \e-\Bigg[ \pd{x}\bar V_{i}(t,x)f_i(t,x) + \pd{t}\bar V_{i}(t,x) \nonumber\\
  & -
  \dot \rho_i(t) + \mu_i(t,x)(\rho_i(t) - \bar V_{i}(t,x))\nonumber \\
  & + \ell_i(t,x)(t-t_i)(t_{i+1}-t) \Bigg] \in \Sigma[t,x],\nonumber \\
  & \ell_i \in \Sigma[t,x]. \nonumber
\end{flalign}
We note that the volume of the set $\fun$ defined by $V(t,x)$
is proportional to:
$$\vol(\fun) \propto \int_{\ti}^{\tf} \sqrt{\frac{\rho(t)^n}{\det(P_0(t))}}
\; dt.$$
As a surrogate for this cost function, we optimize the linear cost:
$$h(\rho) = \int_{\ti}^{\tf} \rho(t)\; dt$$
which can be computed exactly.
% $$h(\rho) = \sum_{i\in \N} \frac{(t_{i+1}-t_i)\rho_i(t_i)$$
% where, of course, a finer sampling that $\{t_i\}_{i\in \N}$ could be
%   used at virtually no computational cost.

We demonstrate that if a given set
$\{\rho_i,\ell_i,\mu_i\}_{i\in \N}$ is feasible, then $\rho(t)$
defines a function $V(t,x)$ satisfying the conditions of Lemma
\ref{eq:posinv}.  First, we see that $\rho$ is constrained to be
continuous and positive so that $V(t,x)$ will be piecewise continuous,
piecewise continuously differentiable and positive.

Writing $V_i = \frac{\bar V_{i}}{\rho_i}$, we have:
\begin{equation}
\frac{\partial}{\partial x}V_i(t,x) f_i(t,x) + \frac{\partial}{\partial
  t}V_i(t,x)
\end{equation}
equivalent to:
\begin{flalign*}
\frac{1}{\rho_i(t)}\Bigg[\frac{\partial}{\partial
    x}\bar V_{i}(t,x) f_i(t,x) +
  \frac{\partial}{\partial
  t}\bar V_{i}(t,x) - \dot\rho_i(t) \frac{\bar V_{i}(t,x)}{\rho_i(t)} \Bigg].
\end{flalign*}
As a result, for $t \in [t_i,t_{i+1}]$ and $x$ such that $\bar V_i(t,x) =
\rho_i(t)$ (equivalently, $x \in \partial\Omega_{t}$)  the optimization
\eqref{eq:bilinear} verifies:
$$\frac{\partial}{\partial x}V_{i}(t,x) f_i(t,x) + \frac{\partial}{\partial
  t}V_i(t,x)  <  \e \rho(t)^{-1}.$$
Note that these constraints verify that:
$$\frac{\partial}{\partial x}V(t,x) f(t,x) + \frac{\partial}{\partial
  t}V(t,x) < 0$$ where we impose the constraint on both left and right
derivatives with respect to time where $\frac{\partial}{\partial t}$
is not continuous.

This optimization is unfortunately bilinear in the coefficients of
$\rho$ and the Lagrange multipliers.  However, this is amenable to an
alternating search.  Given a feasible $\rho(t)$, we can compute
Lagrange multipliers.  Then, holding these multipliers fixed we can
attempt to improve $\rho(t)$.  This is described in detail below.

\subsection{The Multiplier-Step: Finding LaGrange Multipliers}
For a fixed $\rho(t)$ we can compute Lagrange multipliers via the
following set of optimizations. For each interval $[t_i,t_{i+1}]$ we
optimize over two polynomials, $\mu_i$ and $\ell_i$ in $\RR[t,x]$, and
a slack variable $\gamma_i$:
\begin{flalign}\label{eq:fixedrho}
  \minimize_{\gamma_i,\ell_i,\mu_i} \quad & \gamma_i\\
  \subjto \quad & \gamma_i-\Bigg[ \pd{x}\bar V_{i}(t,x)f_i(t,x) + \pd{t}\bar V_{i}(t,x) \nonumber\\
  & -
  \dot \rho_i(t) + \mu_i(t,x)(\rho_i(t) - \bar V_{i}(t,x))  \nonumber \\
  & + \ell_i(t,x)(t-t_i)(t_{i+1}-t)\Bigg] \in \Sigma[t,x], \nonumber\\  
  & \ell_i(t,x) \in \Sigma[t,x].\nonumber
\end{flalign}
These programs can be computed in parallel.  If each $\gamma_i$ is
negative, then the set $\{ \mu_i,\ell_i,\rho_i\}_{i\in \N}$ is a
feasible point for the optimization \eqref{eq:bilinear}. 

An obvious question is how to first obtain a feasible $\rho(t)$.
Motivated by Section \ref{sec:posinv}, we suggest the following search.
We search over a positive constant $c \geq 0$, and take $\rho(t)$ to
be a continuous piecewise polynomial approximation:
$$\rho(t) \approx \exp\left(-c \frac{\tf-t}{\tf-\tau}\right).$$
If a given choice of $c$ does not verify a funnel (i.e. if the
optimal values of \eqref{eq:fixedrho} are not all negative), we iteratively
increase the value of $c$, and potentially the number of knot point of
$\rho(t)$.

% Our optimization procedure involves solving a bilinear Sums-of-Squares
% optimization, which will generally be non-convex. Our approach is
% trajectory based. 

% % \item Find a symmetric positive definite matrix $P_\f
% % \in \R^{n \times n}$ describing an ellipse $\B_\f = \{ x \; | \;
% % \|x-x_G\|_{P_\f}^2 \leq 1\} \subset \goal$.  Maximize the volume of
% % this ellipse.
% % \end{enumerate}

% Given an initial quadratic form $P_0(t): [\ti,\tf] \mapsto \SS_+^n$,
% we construct a more general quadratic form using a piecewise
% polynomial function $\rho: [\ti,\tf] \mapsto (0,\infty)$ by taking:
% $$P(t) = \frac{P_0(t)}{\rho(t)}$$
% With this parameterization we have:
% $$\partial \Omega_1(t) = \{ x \; | \; \| x-\hat x_0\|_{P_0}^2 =
% \rho(t) \}.$$
% The condition:
% $$\pd{x}V(t,x) f(t,x,\pi(t,x)) + \pd{t}V(t,x) < 0$$
% from Lemma~\ref{lem:posinv} is equivalent to:
% \begin{flalign}
%   \rho(t)\left(V_0(t,x) + \pd{t}V_0(t,x)\right) &\nonumber\\
%   - \dot \rho(t) \pd{x}V_0(t,x) f(t,x,\pi(t,x))&< 0. \label{eq:rhoVdot}
% \end{flalign}
% We have simply differentiated, then multiplied through by
% $\rho(t)^2$.  It is important to note that this condition is linear in $\rho(t)$.

\subsection{The $V$-Step: Improving $\rho(t)$.}

Having solved \eqref{eq:fixedrho} for each time interval, we now
attempt to increase the size of the region verified by $V$ by
optimizing to increase $\rho(t)$.  To do so, we pose the following
optimization with $\{\ell_i,\mu_i\}_{i\in\N}$ fixed from a solutions
to \eqref{eq:fixedrho}:
\begin{flalign}\label{eq:fixedlag}
  \maximize_{\{\rho_i\}_{i\in\N}} \quad & h(\rho) \\
  \subjto \quad & \rho_{N-1}(\tf) = 1,\nonumber\\
  &\forall \; i \in \N:\nonumber\\
  & \rho_i(t) \in \Sigma[t], \rho_i(t_{i+1}) = \rho_{i+1}(t_{i+1}),\nonumber\\
  &\e-\Bigg[ \pd{x}\bar V_{i}(t,x)f_i(t,x) + \pd{t}\bar V_{i}(t,x) \nonumber\\
  & -
  \dot \rho_i(t) + \mu_i(t,x)(\rho_i(t) - \bar V_{i}(t,x))\nonumber \\
  & + \ell_i(t,x)(t-t_i)(t_{i+1}-t) \Bigg] \in \Sigma[t,x].\nonumber
\end{flalign}
So long as the class of $\rho_i(t)$ includes the $\rho_i(t)$ used in
the optimizations \eqref{eq:fixedrho}, this optimization will be
feasible, and can only improve the achieved value of $h(\rho)$.

%%% Local Variables:
%%% mode: latex
%%% TeX-master: "ifac2010"
%%% End:

%% file: ifac2010samp.tex
\label{sec:samp}
In terms of computational complexity, the most immediate limiting factor to the above approach is performing Sum-of-Squares optimization is the degree of $t$ in the polynomials being constrained.  We now discuss an approximation to verifying the funnels based on sampling in time.

The proposed approximation validates the conditions of Lemma \ref{lem:posinv} only at finely sampled points.  For each interval $[t_i,t_{i+1}]$ from the above formulation,  we choose a finer sampling $t_i = \tau_{i1} < \tau_{i2} < \ldots < \tau_{iM_i} = t_{i+1}$.
For each such $\tau_{ij}$, we define a Lagrange multiplier $\mu_{ij} \in \RR[x]$.  We pose the bilinear SOS optimization:
\begin{flalign}\label{eq:bilinearsamp}
  \maximize_{\{\rho_i,\mu_{ij}\}_{i\in\N}} \quad & h(\rho)\\
  \subjto \quad  & \rho_{N-1}(\tf) = 1,   \nonumber\\
  &\forall \; i \in \N, \forall\; j \in \{1,\ldots,M_i\}:\nonumber\\
  & \rho_i(t) \in \Sigma[t], \rho_i(t_{i+1}) = \rho_{i+1}(t_{i+1}),\nonumber\\
  & \e- \Bigg[ \pd{x}\bar V_{i}(\tau_{ij},x)f_i(\tau_{ij},x)\nonumber\\
  & + \pd{t}\bar V_{i}(\tau_{ij},x) -  \dot \rho_i(\tau_{ij}) \nonumber\\
  & + \mu_{ij}(x)(\rho_i(\tau_{ij}) - \bar V_{i}(\tau_{ij},x)) \Bigg] \in \Sigma[x]. \nonumber
\end{flalign}
This optimization removes all algebraic dependence on $t$, however it does not in general provide an exact certificate.  One would hope that with sufficiently fine sampling one will recover exactness.  A partial result to this effect exists.

\begin{lem}
  Let $V_i: [t_i,t_{i+1}]\times \RR^n$ and $f_i: [t_i,t_{i+1}] \times \RR^n \mapsto \RR^n$ be continuously differentiable functions of $t,x$.  Further, say that for all $ t \in [t_i,t_{i+1}]$ and $x$ such that $V_i(t,x) = 1$, $\pd{x}V(t,x) \neq 0$.  
  
  Then,if there exists a $t \in [t_i,t_{i+1}]$ and $x \in \RR^n$ such that $V(t,x) = 1$ and:
  $$g(t,x) := \pd{x} V_i(t,x)f_i(t,x) + \pd{t} V_i(t,x) > \d$$
  
  for some positive $\d$, then there exists an open interval around $t$ such that for $\tau$ in that interval there exists a $y \in \RR^n$ with:
  $$V(\tau,y) = 1, \textrm{ and } g(\tau,y) > 0.$$
\end{lem}
\begin{pf}
  As $g$ is continuous in $x$ there exists an $\eta > 0$ such that for all $z$ in  $B(x,\eta) = \{ z \; | \; \|x-z\| < \eta\}$ we have $g(t,z) > \d/2$.  As $\pd{x} V_i(t,x) \neq 0$, there exists a $z_1, z_2 \in B(x,\eta)$ such that $V_i(t,z_1) = 1 + \e_1$ and $V_i(t,z_2) = 1 - \e_2$  for some constants $\e_1,\e_2 > 0$.  As $\pd{t} V_i$ is continuous, and $B(x,\eta)$ is bounded, $\pd{t}V_i$ is bounded.  So there exists an interval around $t$ such that, for any $\tau$ in this interval, $V_i(\tau,z_1) > 1, V_i(\tau,z_2) < 1$. As $g$ is continuous in $t$ there exists a sub-interval such that for $\tau$ in this interval additionally $g(\tau,z) > 0$ for all $z \in B(x,\eta)$.  As $V_i(t,x)$ is continuous, there exists a $\theta \in (0,1)$ such that for $y = \theta z_1 + (1-\theta)z_2$ we have $V_i(\tau,y) = 1$. But, as $B(x,\eta)$ is convex, $y \in B(x,\eta)$ so that $g(\tau,y) > 0$.
\end{pf}
Clearly the function $V_i$ we have considered thus far satisfy the requirements of this Lemma.  The result suggests that given a $V(t,x)$ that does not satisfy the conditions of  Lemma \ref{lem:posinv}, there exists a sufficiently fine uniform sample spacing such that $V(t,x)$ will not be in the feasible set of \eqref{eq:bilinearsamp}.
This result is partial in that it does not construct a sufficient sampling bound from computable quantities, and further applies only to a fixed $V(t,x)$, whereas we optimize over $V$ with fixed sampling intervals.  We are currently studying ways of integrating more constructive bounds into the optimization \eqref{eq:bilinearsamp}.

We use an analogous strategy of bilinear alternation to approach \eqref{eq:bilinearsamp}.  The same strategy is used to find an initial feasible $\rho(t)$, and then Lagrange multipliers and $\rho(t)$ are iteratively improved.

%% file: ifac2010numerics.tex
We first illustrate the general procedure with a one-dimensional
polynomial system.  Our second example is an idealized three degree of
freedom satellite model.  For this second example we compare
numerically the proposed techniques.

\subsection{A One-Dimensional Example}
We examine a one dimensional time-varying polynomial differential equation:
\begin{align}\label{eq:oned}
  \frac{d}{dt} x(t) = f(t,x(t)) = x-\frac{1}{2}x^2 + 2t-2.4t^3,
\end{align}
over the interval $t \in [-1,1]$.  Our goal region is $\goal = [0,1]$.  As the system has no control input, we can solve nearly exact bounds on the region $\fun^\star \subset [-1,1] \times \RR$ which flows into the goal by computing the solutions to \eqref{eq:oned} with final value conditions $x(1) = 1$ and $x(1) = 0$.

To find our inner approximation $\fun \subset \fun^\star$, we compute
an approximate numerical trajectory, $\hat x_0(t)$ with final value
$\hat x_0(1) = x_\f = 0.5$.   We take $P_f = 4$ so that $\E_\f = \{ x
\; | \; |x-x_\f|_{P_\f}^2 \leq 1\}  = [0,1]$.  We numerically solve the Lyapunov equation \eqref{eq:lyap} to determine $P(t)$.

\begin{figure}
  \includegraphics[width=\columnwidth]{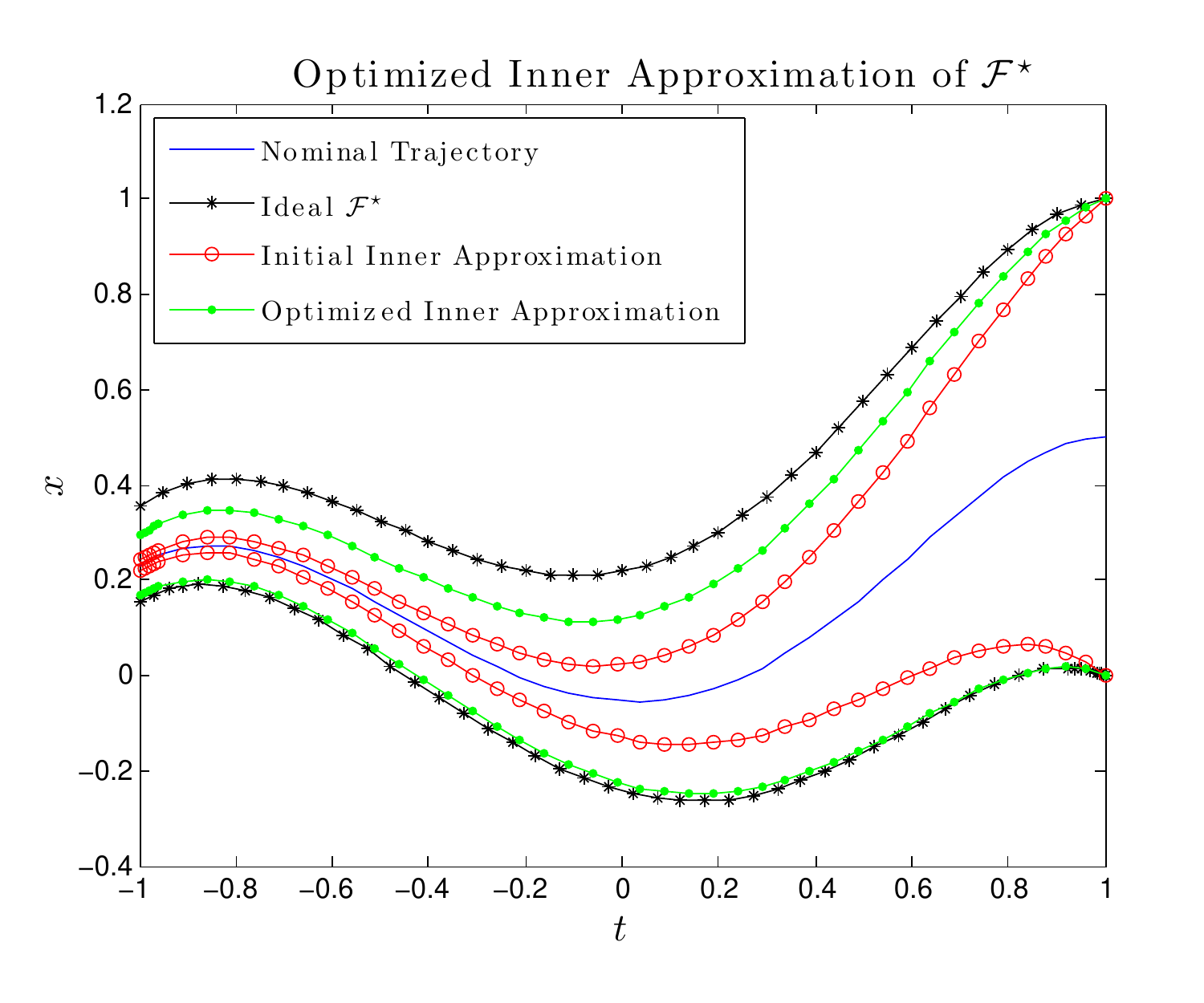}
  \caption{The ideal set $\fun$ and inner approximations calculated by the method. Surrounding the nominal trajectory (solid blue) are time-varying intervals.  An initial candidate Lyapunov function (red open circle) is then improved via the bilinear optimization (solid green circle).  In this case, a single step of alternation provided a certificate tight to the known bounds (black stars).  Note that the certificate is symmetric about the trajectory, and as a result is generally sub-optimal.}
  \label{fig:oned}
\end{figure}

We use $N = 40$ knot points, $\{t_i\}_{i=1}^N$, chosen to be the
steps of a the variable
time-step integration of the differential equation Lyapunov.
We interpolate $\hat x_0(t)$ with a piecewise cubic polynomial and
$P(t)$ with a piecewise linear function.  To find our initial
candidate Lyapunov function, we begin by taking $\rho(t)$ to be a
piecewise linear interpolation of $\exp\left(\frac{c(t-1)}{2}\right)$,
for $c \geq 0$.  Taking $c=4$ provides a feasible candidate Lyapunov
function.  This feasible solution is then improved by bilinear
alternation.  Both the initial and optimized sets are plotted against
the known bounds in Figure \ref{fig:oned}.

After a single bilinear alternation, a tight region is found.  Note that the symmetry of the Lyapunov function around the trajectory restricts the region being verified.  Additional trajectories could be used to continue to grow the verified region.

\subsection{Trajectory Stabilization of Satellite Dynamics}
We next evaluate a time-varying positively invariant region around a feedback stabilized nominal trajectory.  In past work, \cite{Tedrake10}, it was demonstrated how trajectory optimization and randomized search can be combined with such certificates to approximate the controllable region for a smooth nonlinear system. 

We examine the stabilization of a nominal trajectory for a rigid body floating in space subject to commanded torques. The state of the system is given in terms of the angular velocity of the body about is principal axes, $\omega \in \RR^3$, and the Modified Rodriguez parameters, $\sigma \in \RR^3$.  Any trajectory which excludes full rotations can be represented by this projection of the quaternion representation of orientation.  The kinematic equations are given by:
\begin{equation}\label{eq:satkin}
  \dot \sigma = \frac{1}{4} \left((1-\|\sigma\|^2)I + 2\sigma\sigma'
    - 2\begin{bmatrix} 0 & \sigma_3 & \sigma_2 \\ \sigma_3 & 0 & \sigma_1 \\ \sigma_2 & \sigma_3 & 0
    \end{bmatrix}\right)\omega.
\end{equation}
The dynamic equations are given by:
\begin{equation}\label{eq:satdyn}
  H \dot \omega = -(\omega \times H\omega) + u
\end{equation}
where $H = H' > 0$ is the diagonal, positive definite inertia matrix
of the system and $u \in \RR^3$ is a vector of torques.  In our
example $H = \textrm{diag}([5,3,2])$. The state of the system is $x = [ \sigma' , \omega']' \in \RR^6$.  Together, \eqref{eq:satkin} and \eqref{eq:satdyn} define controlled dynamics:
\begin{equation} \label{eq:ctldyn}
  \dot x(t) = f_0(x(t),u(t)).
\end{equation}
We design a control policy $u(t) = \pi(t,x(t))$ such that the closed loop system,
\begin{equation}
  \dot x(t) = f(t,x(t)) = f_0(x(t),\pi(t,x(t))),
\end{equation}
satisfies the assumptions of our method.

Our goal region is defined by an ellipse centered on the origin, described by the positive definite matrix:
$$P_G = \begin{bmatrix}
   36.1704   &0    &     0 &  12.1205 & 0 & 0 \\
   0&   17.4283  &       0 &  0 &  7.2723 &   0 \\
   0&       0 &   9.8911  &       0 &        0 &   4.8482\\
   12.1205 &  0 &0    &9.1505&    0 &        0\\
   0&    7.2723  &0    &0   & 7.3484 &        0\\
   0 &    0          &4.8482&         0 &        0&    6.2557
 \end{bmatrix}.$$ 
 The ellipsoidal region, $\goal = \{ x \; | \; \|x\|_{P_G}^2 \leq 1\}$ was computed numerically as an inner-approximation of the region of attraction for the system after feedback stabilization of the origin.  

We begin with a nominal command:
$$ u_0(t) = \frac{1}{100}t(t-5)(t+5)\begin{bmatrix}-1\\-1\\1\end{bmatrix}.$$
on the interval $t \in [0,5]$.  We compute the solution, $x_0(t)$ to \eqref{eq:ctldyn} with $u(t) \equiv u_0(t)$ and $x(5) = 0$.  Next, we design a time-varying LQR controller around the trajectory based on the dynamics linearized about the trajectory.
Given two cost matrices, $R \in \SS_+^m$ and $Q \in \SS_+^n$, we solve the Riccati differential equation:
$$ -\dot S^\star = A'S^\star + S^\star A + Q - S^\star BR^{-1}B'S^\star, \quad S^\star(5) = P_\f$$
where $P_f = 1.01 P_\G$ and:
$$A(t) = \frac{\partial}{\partial x}f_0(x_0(t),u_0(t)), \quad B(t) = \frac{\partial}{\partial u} f_0(x_0(t),u_0(t)).$$

This procedure gives us a time-varying gain matrix:
$$K(t) = R^{-1}B(t)'S^\star(t)$$
The ideal policy is given by:
$$\pi^*(t,x(t)) = u_0(t) -  K(t)(x- x_0(t)).$$
To force $\pi(t,x)$ to be piecewise polynomial in $t$ and polynomial in $x$ we take a piecewise constant approximation of $\hat K$ of $K$ and a piecewise cubic approximation $\hat x_0(t)$ of $x_0(t)$.  Our control policy is then:
$$\pi(t,x(t)) = u_0(t) -  \hat K(t)(x- \hat x_0(t)).$$

We now examine the closed loop dynamics using the Lyapunov function:
$$V(t,x) = \| x - \hat x_0(t)\|_{S_0(t)/\rho(t)}^2$$
where $S_0(t)$ is a piecewise linear approximation of $S^\star(t)$.
All of the above approximations were piecewise with $N = 49$ knot points chosen by a variable time-step integration of the Riccati differential equation.

\begin{figure}
  \centering
  \includegraphics[width=\columnwidth]{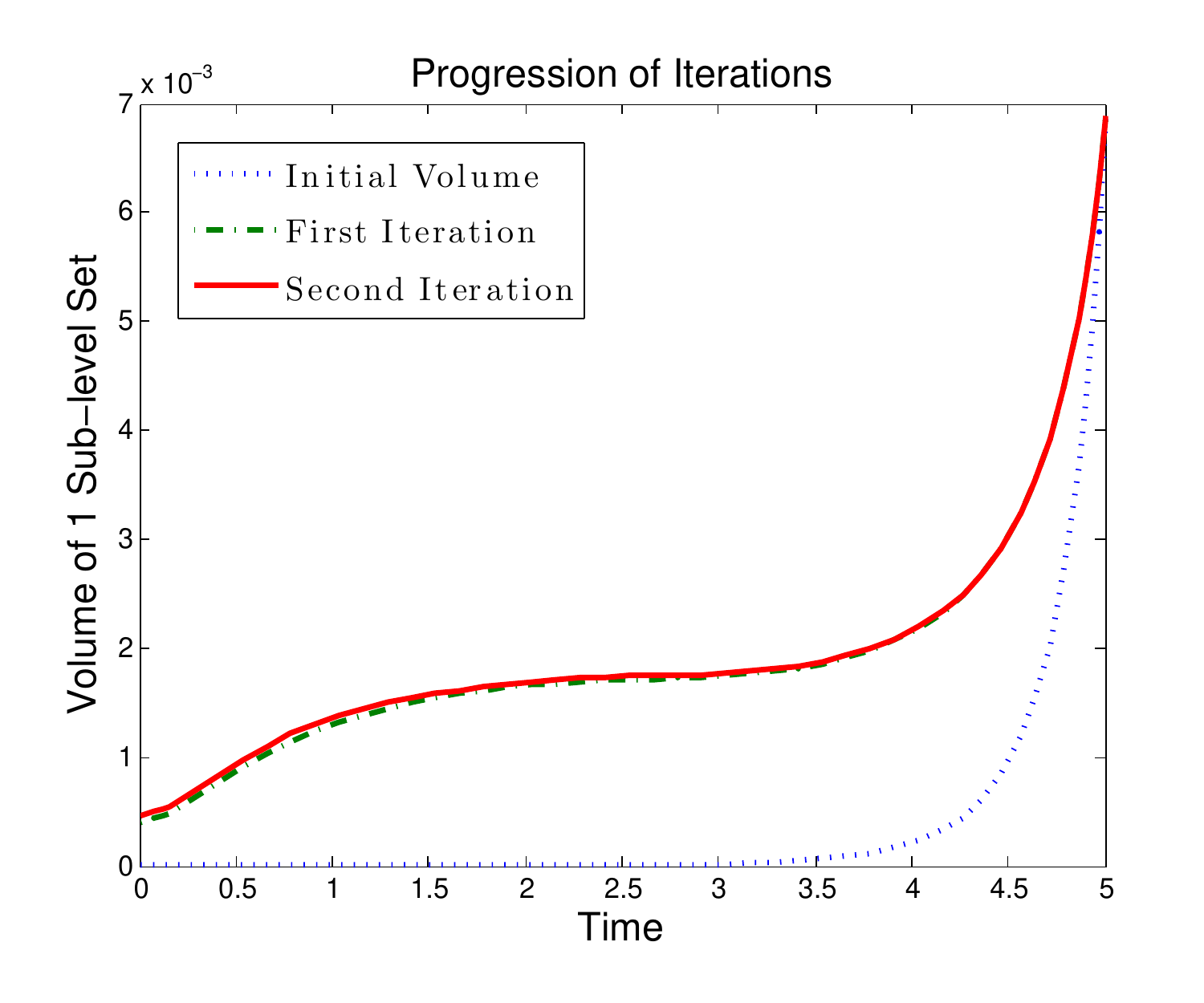}
  \caption{The result of successive iterations optimizing the volume of verified funnels.  After two iterations, the methods, the volume ceases to improve.}
  \label{fig:iterations}
\end{figure}

\begin{figure}
  \centering
  \includegraphics[width=\columnwidth]{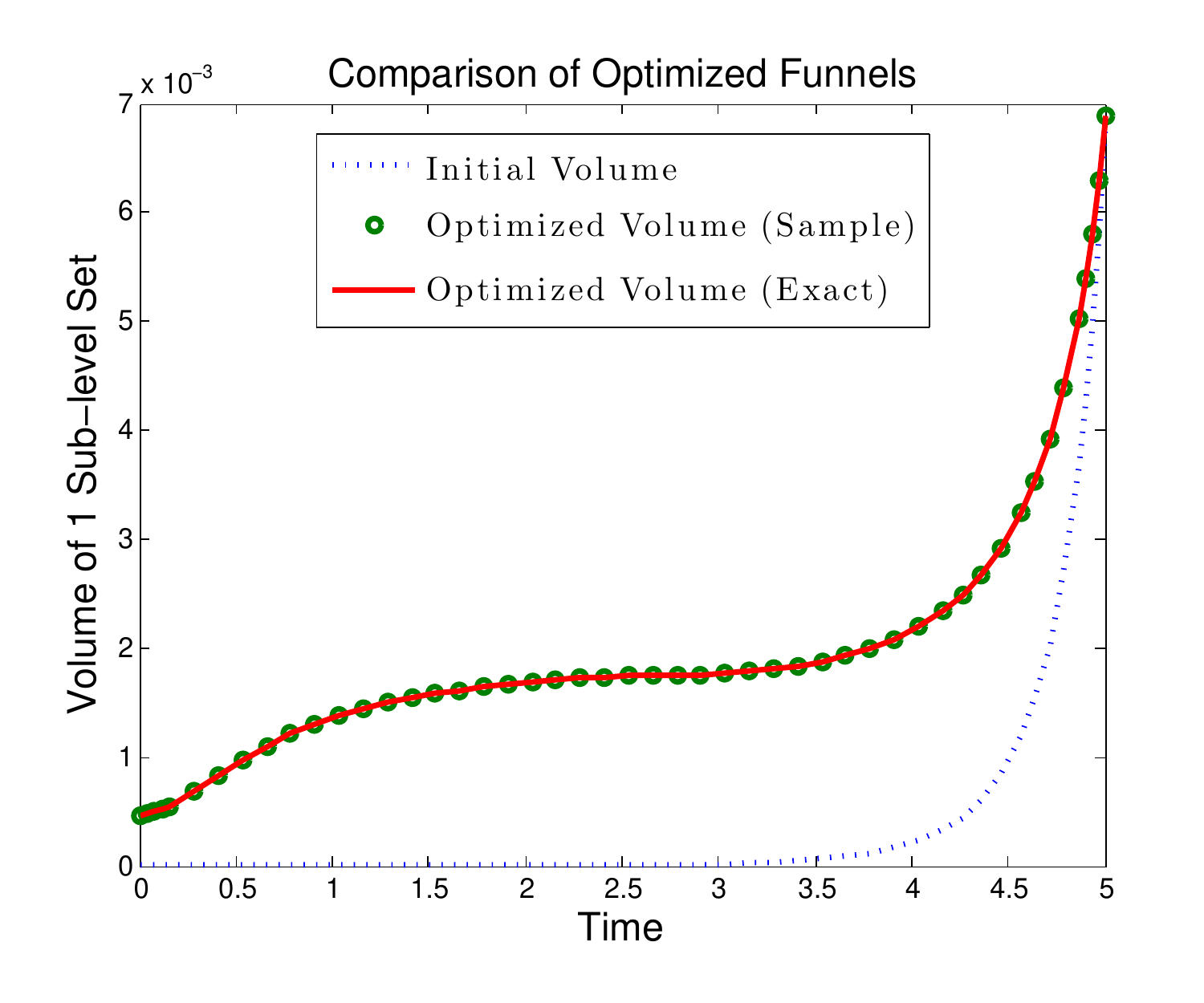}
  \caption{Comparison of optimized volumes using the exact method and time-sampled method.}
  \label{fig:comparison}
\end{figure}

We compare the time required to compute an exact certificate versus
the sample based approximation verifying necessary conditions at $M_i
= 10$ points in each interval $[t_i,t_{i+1}]$.   In the case of
computing exact certificates, the expression for $\frac{d}{dt} V(t,x(t))$ is degree 7 in $t$. For both processes, we found that the exponentially decaying initial $\rho(t)$ with $c = 3$ was feasible.  
Figure \ref{fig:iterations} shows the progress of two iterations using
the exact method against the initial $\rho(t)$.  The figure plots the
volume of $\Omega_t$ for each $t \in [0,5]$.
We have found that typically very few iterations are required for the
procedure to converge; indeed, in this case the first iteration was
already very good.
Figure \ref{fig:comparison} compares the results from two iterations of the time-sampled and exact methods.  In this instance they are nearly identical.
Table \ref{tab:speed} compares the run times. The Multiplier Step consisted of 48 independent SDPs for the exact method and 480 independent SDPs for the sample based method.  These computations were performed on a 12-core 3.33 GHz Intel Xeon computer with 24 Gb of RAM.
These Multiplier Step programs could be trivially parallelized for speedup. 

\begin{table}
  \centering
  \caption{Runtime comparison for exact and sample based approaches.  }
  \label{tab:speed}
  \begin{tabular}{|l|rr|rr|}
    \hline
    & \multicolumn{2}{c}{Multiplier Step (sec)} & \multicolumn{2}{|c|}{$V$ Step (sec)}\\
    \hline
    Sampled & 111& 112 & 220&220\\
    Exact& 5316 & 5357 &1336 & 1420\\
    \hline
  \end{tabular}
\end{table}

\begin{figure}
  
\end{figure}

%%% Local Variables:
%%% mode: latex
%%% TeX-master: "ifac2010"
%%% End:

%% file: ifac2010discussion.tex
We now discuss simple extensions and implementation details for the method.

% The simplest such extension is search over $P(t)$ directly as a time-varying positive-definite matrix function.  Similarly, one could search over functions of the form:
% $$V(t,x) = \begin{bmatrix}1 \\ x-\hat x_0(t)\end{bmatrix}'
% M(t)
% \begin{bmatrix}1 \\ x-\hat x_0(t)\end{bmatrix}, \\
% M(t) = \begin{bmatrix} c(t) & f(t)' \\ f(t) & P(t)\end{bmatrix} \geq 0$$
\subsection{Sampling in Time without Splines}
When applying the time-sampling method described in this paper, one is no longer bound to finding a polynomial Lyapunov function.  As a result, one case use the numerically determined $x_0(t)$ and $P_0^\star(t)$ instead of interpolations.  We make use of the splined quantities in both the sampled based and exact methods in this work only for the sake of comparison.
\subsection{Verifying Families of Funnels}
For several applications, it may prove useful to verify that:
$$\pd{x} V(t,x)f(t,x) + \pd{t} V(t,x) < 0$$
for all level-sets in a range $V(t,x) \in [a,1]$, with $0 < a < 1$.  For example: when performing real-time composition of preplanned trajectories, it may not be known in advance what the goal region is for each trajectory segment. By verifying funnels for a whole family of goal regions, the best verified funnel can be chosen at execution time.
\subsection{More General Lyapunov Functions}
In this paper work we have restricted ourselves to  rescalings of a time-varying quadratic form, centered near trajectories.  Our first numerical example demonstrated how this class can be quite conservative.  In principle more general polynomial Lyapunov function can be addressed with the method presented.  Using a richer class of polynomial Lyapunov functions has proven advantageous for time-invariant region-of-attraction analysis (for example, see \cite{Topcu08}).  We are currently investigating extending the given algorithm for this case.
\subsection{Stability of Limit Cycles}
The authors have also adapted the procedure proposed in this paper to verify regions of attraction for limit cycles of hybrid systems (see \cite{Manchester10a, Manchester10b}).

%%% Local Variables:
%%% mode: latex
%%% TeX-master: "ifac2010"
%%% End: